\RequirePackage{ifpdf}
\ifpdf 
\documentclass[pdftex]{sigma}
\else
\documentclass{sigma}
\fi

\begin{document}

\allowdisplaybreaks

\renewcommand{\PaperNumber}{118}

\FirstPageHeading

\renewcommand{\thefootnote}{$\star$}

\ShortArticleName{On Gauss--Bonnet Curvatures}

\ArticleName{On Gauss--Bonnet Curvatures\footnote{This paper is a
contribution to the Proceedings of the 2007 Midwest
Geometry Conference in honor of Thomas~P.\ Branson. The full collection is available at
\href{http://www.emis.de/journals/SIGMA/MGC2007.html}{http://www.emis.de/journals/SIGMA/MGC2007.html}}}

\Author{Mohammed Larbi LABBI}

\AuthorNameForHeading{M.L. Labbi}

\Address{Mathematics Department, College of Science, University of Bahrain, 32038 Bahrain}
\Email{\href{mailto:labbi@sci.uob.bh}{labbi@sci.uob.bh}}

\ArticleDates{Received August 27, 2007, in f\/inal form November
15, 2007; Published online December 11, 2007}

\Abstract{The $(2k)$-th Gauss--Bonnet curvature is a generalization to higher dimensions
of the $(2k)$-dimensional Gauss--Bonnet integrand, it coincides with the usual scalar curvature for $k=1$.
The Gauss--Bonnet curvatures are used in theoretical
physics  to describe gravity in higher dimensional
space times where they are known as the Lagrangian of Lovelock gravity, Gauss--Bonnet Gravity and Lanczos gravity.
In this paper we present various aspects of these curvature invariants and  review  their
 variational properties. In particular, we discuss  natural generalizations of the Yamabe
problem, Einstein metrics and minimal submanifolds.}


\Classification{53C20; 53C25}

\section[An introduction to Gauss-Bonnet curvatures]{An introduction to Gauss--Bonnet curvatures}

We shall present in this section several approaches to the Gauss--Bonnet curvatures. For precise def\/initions and examples the reader is encouraged to
 consult \cite{Labbidoubleforms,LabbiPacific,LabbiHabilitation}.

\subsection[Gauss-Bonnet curvatures vs.  curvature invariants of Weyl's tube formula]{Gauss--Bonnet curvatures vs.\  curvature invariants\\ of Weyl's tube formula}

In a celebrated paper \cite{Weyl} published in 1939, Hermann Weyl proved that the volume of a tube of radius $r$ around an
embedded compact $p$-submanifold $M$ of the
$n$-dimensional  Euclidean space is a polynomial in the radius of the tube as follows:
\begin{gather*}
{\rm Vol}({\rm tube}(r))= \sum_{i=0}^{[p/2]} C(n,p,i) H_{2i}r^{2i},
\end{gather*}
where $C(n,p,i)$ are constants which only depend on the dimension  and the codimension of the submanifold $M$, and $H_{2i}$ are integrals of intrinsic
 scalar curvatures of the submanifold (Gauss--Bonnet curvatures). A good and complete reference about Weyl's tube  formula and related topics is the book of
A. Gray \cite{Graybook}.

It turns out that
$H_0$ is the volume of the submanifold,
$H_2$ is the integral of the usual scalar curvature of the submanifold,
the integrand in $H_4$ is quadratic in the Riemann tensor and was introduced by Lanczos in 1932 as a possible substitute to
Hilbert's Lagrangian in general relativity. The top $H_p$ is up to a constant the Euler--Poincar\'e characteristic of the submanifold if $p$ is even.

All the (total Gauss--Bonnet) curvatures $H_{2i}$ have important applications in theoretical physics, particularly in (brane world) cosmology.
They are by nowadays the subject of intensive studies, where they are  known as the Lagrangian of Lovelock gravities or Gauss--Bonnet gravities,
see \cite{Chijap,Deruelle-Madore,Madore} and the references therein.

\subsection[Gauss-Bonnet curvatures vs. Gauss-Bonnet integrands]{Gauss--Bonnet curvatures vs.\ Gauss--Bonnet integrands}

\subsubsection{From Gaussian curvature to the scalar curvature}

Recall that for a compact $2$-dimensional Riemannian manifold $(M,g)$ (a surface) the classical Gauss--Bonnet formula states that the Euler-Poincar\'e characteristic
of $M$ (which is a topological invariant) is determined by the geometry of $(M,g)$ as an integral of the
 Gaussian curvature  of the metric: the 2-dimensional Gauss--Bonnet integrand.
It is a scalar function def\/ined on the surface and can be naturally generalized to higher dimensional Riemannian manifolds in the following way:

Let $(M,g)$ be a Riemannian manifold of dimension $n\geq 2$.
For  $m\in M$ and for a tangent $2$-plane $P$ to $M$ at $m$ we def\/ine $K(P)$, the sectional
curvature at $P$, to be the Gaussian curvature at $m$ of the surface $\exp_m(V)$,
where ${\rm exp}_m$ is the exponential map and $V$ is a small neighborhood of $0$ in $V$. Recall
that the so-obtained surface is totally geodesic at $m$. In this way, we obtain a function $K$
def\/ined on the $2$-Grassmannian bundle over  $M$. The function $K$ determines   a~unique  symmetric $(2,2)$-double form $R$ that satisf\/ies the f\/irst
Bianchi identity and having $K$ as its
sectional curvature, that is nothing but the standard Riemann curvature tensor. Recall that a symmetric $(2,2)$-double form is a $(0,4)$
tensor  which is skew symmetric in the f\/irst two arguments and in the last two, and that it is symmetric with respect
 to the interchange of the f\/irst two variables with the last two.

Then one can def\/ine the scalar curvature of $M$ by taking the Ricci contraction of $R$ twice. In this sense one can say that
the usual scalar curvature is a natural  generalization of the two dimensional Gauss--Bonnet integrand to higher dimensions.

\subsubsection[From higher Gauss-Bonnet integrands to Gauss-Bonnet curvatures]{From higher Gauss--Bonnet integrands to Gauss--Bonnet curvatures}

For a compact $(2p)$-dimensional Riemannian manifold $(M,g)$ the generalized Gauss--Bonnet theorem states that the Euler--Poincar\'e characteristic
of $M$ (which is a topological invariant) is determined by the geometry of $(M,g)$ as an integral of a certain curvature of the metric:
\[\chi(M)=c(p)\int_M h_{2p}{\rm \, dvol},\]
where $c(p)$ is a constant and $h_{2p}$ is a scalar function on the manifold def\/ined using the Riemann curvature tensor of $(M,g)$:
the $(2p)$-th Gauss--Bonnet integrand, called also the Lipschitz--Killing curvature~\cite{Fenchel}.

Using the same idea as above, we generalize the $(2p)$-th Gauss--Bonnet integrand (or the $(2p)$-th Killing--Lipschitz curvature) to
dimensions higher than $(2p)$  as follows.

Let $(M,g)$ be a Riemannian manifold of dimension $n\geq 2p$.
 For  $m\in M$ and for a tangent $(2p)$-plane $P$ to $M$ at $m$ we def\/ine $K_{2p}$,
 called Thorpe's $(2p)$-th sectional
curvature at $P$ \cite{Stehney}, to be the Gauss--Bonnet integrand at $m$ of the $(2p)$-dimensional
submanifold  $\exp_m(V)$,
where $\exp_m$ and $V$ are as above. Thorpe's tensor $R_{2p}$ of order $(2p)$ is then  def\/ined to be  the unique
symmetric $(p,p)$-double form that satisf\/ies the f\/irst Bianchi identity and with sectional curvature~$K_{2p}$.

Then one can get scalar curvatures (Gauss--Bonnet curvatures) after taking $(2p)$-times the Ricci contraction of $R_{2p}$.

The tensors $R_{2p}$ are determined by the Riemann curvature tensor $R$ in the following way \cite{Thorpe}:
For $u_i$, $v_j$ tangent vectors at $m\in M$, we have
\begin{gather*}
\frac{(2p)!}{2^p} R_{2p}(u_1,\ldots,u_{2p},v_1,\ldots,v_{2p})\\
\qquad{}=\sum_{\alpha, \beta \in S_{2p}}
\epsilon(\alpha)\epsilon(\beta)
 R(u_{\alpha(1)},u_{\alpha(2)},
v_{\beta(1)},v_{\beta(2)}) \cdots  R(u_{\alpha(2p-1)},u_{\alpha(2p)},
v_{\beta(2p-1)},v_{\beta(2p)}).
\end{gather*}
This complicated expression can be considerably simplif\/ied using the exterior product of double forms, see the following subsection.

\subsection{Double forms}

A $(p,q)$-double form
$\omega(x_1,x_2,\ldots,x_p,y_1,y_2,\ldots,y_q)$ on $M$ is at each point of the manifold  a~multilinear form that is skew symmetric with respect
to the interchange of any two among  the f\/irst $p$-arguments (tangent vectors) or the last~$q$. If $p=q$ and
$\omega$ is invariant with the respect
to the interchange of the f\/irst $p$-variables with the last~$p$, we say that $\omega$ is a symmetric $(p,p)$-double form.

For example, the covariant Riemann curvature tensor is a symmetric $(2,2)$-double form, and
Thorpe's tensor $R_{2p}$ is  a symmetric $(2p,2p)$-double form.

Double forms were introduced f\/irst by De Rham and then developed in the seventies of the last century by Thorpe \cite{Thorpe}, Kulkarni \cite{Kulkarni},
Nomizu \cite{Nomizu}, Gray \cite{Gray}, Kowalski \cite{Kowalski}, Nasu \cite{Nasu} etc.

Double forms and the more general multi-forms are recently
studied in theoretical physics, see~\cite{Senovilla} and the references therein.

\subsubsection{Algebraic operations on double forms}

A $(p,q)$-double form can be seen alternatively as a section of the tensor product of the bundle of $p$-forms with the one of $q$-forms.
\begin{enumerate}\itemsep=0pt
\item  The exterior  product  of double forms
 is the natural generalization to double forms of the usual exterior product of dif\/ferential forms:
\[(\theta_1\otimes \theta_2)\cdot (\theta_3\otimes \theta_4)=
    (\theta_1\wedge \theta_3 )\otimes(\theta_2\wedge \theta_4),
    \]
where we denoted the exterior product of double forms following Kulkarni \cite{Kulkarni} by a dot.
 In the following the dot shall be omitted whenever possible. The well known Kulkarni--Nomizu product of symmetric bilinear forms is a special case
of the exterior product, see \cite[p.~47]{Besse}.

\item The generalized Hodge star operator \cite{Labbidoubleforms}
is the natural extension to double forms of the usual Hodge star operator on dif\/ferential forms:
\begin{gather*}
*(\theta_1\otimes \theta_2)=(*\theta_1)\otimes (*\theta_2).
\end{gather*}
Senovilla \cite{Senovilla} considered other interesting extensions of the Hodge star operator to double\-forms and multi-forms by keeping it acting on dif\/ferent
 factors of the tensor product.
\item The inner product of double forms \cite{Labbidoubleforms} is def\/ined by declaring
\begin{gather*}
    \langle \theta_1\otimes \theta_2, \theta_3\otimes \theta_4 \rangle=
    \langle\theta_1, \theta_3 \rangle  \langle\theta_2, \theta_4\rangle.
\end{gather*}
\end{enumerate}

The exterior product of double forms has the advantage that it makes easier many complicated expressions of
Riemannian geometry, as illustrated by the examples below.
 In the following  let
$(M,g)$ denote a Riemannian manifold of dimension $n$ and  $R$ its Riemann curvature tensor seen as a $(2,2)$-double form.
\begin{itemize}\itemsep=0pt
\item Thorpe tensors are just given by \cite{Kulkarni, Nasu}
\begin{gather*}
R_{2p}=\frac{2^{p}}{(2p)!}R^{p},
\end{gather*}
where of course $R^p$ is the exterior product of the Riemann curvature tensor $R$ seen as a $(2,2)$-double form. In particular,
$R^{n/2}$ determines the Gauss--Bonnet integrand if the dimension $n$ of the manifold is even.
\item
The $(2k)$-th Gauss--Bonnet curvature of $(M,g)$ can be alternatively def\/ined by \cite{Labbidoubleforms}:
\begin{gather}\label{GaussBonneteven}
h_{2p}= *\frac{1}{(n-2p)!}(g^{n-2p}R^p).
\end{gather}

\item The curvature operator of the classical  Weitzenb\"{o}ck formula acting on $p$-forms is determined
 by the following double form \cite{Bourguignon, LabbiWeitzenbock}.
\[\left\{\frac{g{\rm Ric}}{(p-1)}-2R\right\}\frac{g^{p-2}}{(p-2)!},\]
where ${\rm Ric}$ denotes the Ricci curvature of $(M,g)$.
\item The exterior product $g^k=g\cdots g$ determines the canonical inner product of
 dif\/ferential $k$-forms over $M$.

\item  Gauss equation for a hypersurface of the Euclidean space which relates the curvature tensor~$R$ of the hypersurface
to its second fundamental form $B$ can be just written as \cite{Graybook}
\[ R=\frac{1}{2}B^2.\]

\end{itemize}

\subsection[Gauss-Bonnet curvatures vs. symmetric functions in the eigenvalues of the shape operator of a
 hypersurface of the Euclidean space]{Gauss--Bonnet curvatures vs.\ symmetric functions in the eigenvalues\\ of the shape operator of a
 hypersurface of the Euclidean space}

Let $g$ and $B$ denote respectively the f\/irst and  second fundamental forms of a hypersurface of the Euclidean space.
The symmetric functions in the eigenvalues of the operator corresponding to $B$ can be nicely re-formulated
using the exterior product and the generalized Hodge star operator as follows \cite{Labbi2kminimal}
\[
s_k=\frac{1}{k!(n-k)!}*(g^{n-k}B^k).
\]
In particular, if $k=2p$ is even,
Gauss equation shows that $R=\frac{1}{2}B^2$.
Therefore, all the even powers of $B$ are then intrinsic and consequently $s_{2p}$ is also intrinsic and coincides up to a constant with the Gauss--Bonnet
curvature of the hypersurface as follows:
\[ s_{2p}=\frac{2^p}{(2p)!(n-2p)!}*(g^{n-2p}R^p)=\frac{2^k}{(2k)!}h_{2k}.\]
Note that if $k=2p+1$ is odd then $s_{2p+1}$ is  no longer intrinsic:
\[ s_{2p+1}=*\frac{g^{n-2p-1}B^{2p+1}}{(n-2p-1)!(2p+1)!}=
*\frac{2^pg^{n-2p-1}R^p B}{(n-2p-1)!(2p+1)!}.\]
The previous formula allows one to def\/ine the Gauss--Bonnet curvatures of odd order for an~arbitrary submanifold as follows:
\begin{definition}[\cite{Labbi2kminimal}]
Let $(M,g)$ be an arbitrary $n$-submanifold of a Riemannian manifold $(\tilde M, \tilde g)$ and $N$ a normal vector to $M$. For
$1\leq 2p+1\leq n$, we def\/ine
the $(2p+1)$-th Gauss--Bonnet curvature of  $(M,g)$ at $N$ by
\begin{gather}\label{GaussBonnetodd}
h_{2p+1}(N)=*\left(\frac{g^{n-2p-1}}{(n-2p-1)!}R^pB_N\right),
\end{gather}
where $B$ denotes the vector valued second fundamental form of $M$,
$B_N(u,v)=\tilde g(B(u,v),N)$ and $R$ is the Riemann curvature tensor of $(M,g)$.
\end{definition}

The $(2p+1)$-th Gauss--Bonnet curvature is a generalization of the usual mean curvature as for  $p=0$ we recover the trace of $B$:
\[ h_1(N)=*\left(\frac{g^{n-1}}{(n-1)!}B_N\right)=cB_N.\]
Furthermore, for a submanifold of the Euclidean space, $h_{2p+1}$ coincides with the higher
$(2p+1)$-th mean curvature def\/ined by Reilly~\cite{Reillyproc}.

\subsection[Other aspects of Gauss-Bonnet curvatures and terminology]{Other aspects of Gauss--Bonnet curvatures and terminology}

The Gauss--Bonnet curvatures, or more precisely closely related invariants, appear naturally in the mathematics and physics literature
under dif\/ferent names.
For instance, in the context of convex sets they (more precisely their integrals) are known as Minkowski's Quermassintegrale,
intrinsic volumes or Steiner functionals.
They are called Lipschitz--Killing curvatures (measures) in the case of sets of positive reach, piecewise linear spaces and subanalytic sets.
A good compte-rendu of these dif\/ferent aspects is the review paper of Bernig~\cite{Berningreview}. These invariants appear also in the study of
the ``expected Euler characteristic'' in the probability literature~\cite{Adler} and in extensions of the Gauss--Bonnet theorem \cite{Albin}.

In theoretical physics Gauss--Bonnet curvatures are known as a Lagrangian of pure Lovelock gravity, Gauss--Bonnet gravity or Lanczos gravity.
Precisely, in dimensions higher
than four the Lagrangian of Lovelock gravity is \cite{Chijap}
\begin{gather}\label{gravity}
L=\sum_{k=0}^{m}c_{2k}h_{2k},
\end{gather}
where $h_0=1$, $c_0$ is the cosmological constant, the higher $c_k$ are arbitrary constants and $h_{2k}$ (resp.~$m$) is the $(2k)$-th
Gauss--Bonnet curvature (resp.\ the integer part of one half the dimension of) the (pseudo) Riemannian
manifold under consideration.

Physicists consider the Lovelock gravity as a fascinating extension of general relativity for dimensions higher than four, see for instance \cite{Chijap}
and the references therein. The Lagrangian of the so-called pure Lovelock gravity is the Gauss--Bonnet curvature
 (that is just one term in the summation in (\ref{gravity})).
The Lagrangian of the Gauss--Bonnet gravity (called also Lanczos gravity) is extensively studied in theoretical physics,
it is obtained by taking only the f\/irst three terms in the summation (\ref{gravity}),
see \cite{Madore,Chijap} and the references therein.

Finally we make the following comments about terminology.
Some authors name all the previous invariants as Killing--Lipschitz curvatures in all categories. In the smooth Riemannian case Killing--Lipschitz
curvatures refer originally to the Gauss--Bonnet integrands \cite{Fenchel}. On the other hand, there is no standard terminology for the scalar
curvatures (that are the Gauss--Bonnet curvatures under study here) generalizing the previous integrands. For instance, in
\cite{Berningreview,Cheeger,Lafontaine} they are still called Lipschitz--Killing curvatures,
in \cite{Berger,Japweyl} they are H.~Weyl's  invariants, in \cite{Graybook, Reillyproc} they are called higher mean curvatures, physicists call them
dimensionally continued (or extended) Euler densities or Gauss--Bonnet integrands. It was a subtle suggestion made to us by an anonymous referee of the
{\it Pacific J. Math.}  to rather call these higher order scalar curvatures in the smooth (Riemannian) case as the Gauss--Bonnet or Gauss--Bonnet--Weyl curvatures.

\section[Einstein-Lovelock tensors]{Einstein--Lovelock tensors}

The usual Ricci curvature tensor $cR$ is the f\/irst Ricci-contraction of the Riemann curvature tensor
$R$. The Einstein tensor is the simplest linear combination of the the Ricci tensor and
the metric tensor to be divergence free, that is $\frac{1}{2}c^2R g-cR$. It is the gradient
of the total scalar curvature seen as a functional on the space of all Riemannian metrics on the manifold under consideration.

In a similar way, we def\/ine a generalized Ricci curvature tensor $c^{2p-1}R^p$ of order $(2p)$
 to be the $(2p-1)$-th Ricci contraction of Thorpe's tensor $R^p$. The Einstein--Lovelock
tensor $T_{2p}$ is a~linear combination of  the $(2p )$-th Ricci tensor
$c^{2p-1}R^p$  and
the metric tensor that is  divergence free. Precisely, we def\/ine the Einstein--Lovelock tensor $T_{2p}$
of order $2p$ by
\begin{gather}\label{Einstein-Lovelock}
 T_{2p}=h_{2p}g-\frac{1}{(2p-1)!}c^{2p-1}R^p.
 \end{gather}

For $p=0$, we set $T_0=g$. For $p=1$, $T_2$ coincides with the usual Einstein tensor. Furthermore, the
  tensor $T_{2k}$ is the gradient of the total $(2k)$-th Gauss--Bonnet curvature
seen as a functional on the space of all Riemannian metrics on
 a given compact manifold, see the next section.

David Lovelock \cite{Lovelock}  proved that any divergence-free symmetric $(0,2)$ tensor built from the metric and its f\/irst two
covariant derivatives are linear
combinations of the tensors $T_{2k}$.

\section[A variational property of the Gauss-Bonnet curvatures]{A variational property of the Gauss--Bonnet curvatures}

On a compact manifold, we have the classical total scalar curvature functional:
\[
S(g)=\int_M {\rm scal}(g) \mu_g.
\]
The gradient of this Riemannian functional is the Einstein tensor: $\frac{1}{2}{\rm scal}g-{\rm Ric}.$
The critical metrics  of $S$ once restricted to metrics with unit volume,
are the
Einstein metrics.

Similar properties hold for the total Gauss--Bonnet curvature functional:
\[ H_{2k}(g)=\int_Mh_{2k}(g)\mu_g,\]
as shown by the following theorem:
\begin{theorem}
  Let $(M,g)$ be a compact Riemannian manifold of  dimension $n$. For each   $k$, such that
 $2\leq 2k\leq n$,
the functional  $H_{2k}$ is differentiable, and at   $g$ we have
  \[ H_{2k}'h=\frac{1}{2}\langle h_{2k}g-\frac{1}{(2k-1)!}c^{2k-1}R^k,h\rangle.\]
In particular, the  gradient of $2H_{2k}$ is  $T_{2k}=h_{2k}g-\frac{1}{(2k-1)!}c^{2k-1}R^k$.
\end{theorem}
The previous theorem were f\/irst proved by Lovelock \cite{Lovelock} using classical tensor analysis, Berger~\cite{Berger} in the case $k=2$,
Patterson \cite{Patterson} as a special case of a more general variational formula due to Muto~\cite{Muto} and Bernig~\cite{Berningvariation} proved
a more general variational formula for the Lipschitz--Killing curvatures of subanalytic sets. The  proof of \cite{Labbivariation} sketched
below is coordinate free and it uses the
the formalism of double forms.

\begin{proof}
We sketch the proof of the theorem, for more details see \cite{Labbivariation}.
First, the directional derivative of the Riemann curvature tensor $R$, seen as
a symmetric double form has the form:
\begin{gather*}
R'h={\rm Exact\,  double \, form}+ {\rm a\, linear\, term\, in\, R},
\end{gather*}
precisely,
\begin{gather*}
R'h=-\frac{1}{4}(D\tilde D+\tilde D D)(h)+\frac{1}{4}F_h(R).
\end{gather*}
Next, we derive the directional derivative of the Gauss--Bonnet curvature  $h_{2k}$ at $g$:
\[
h_{2k}'h=-\frac{1}{2}\langle\frac{c^{2k-1}}{(2k-1)!}R^k,h\rangle-
  \frac{k}{4}(\delta\tilde\delta +\tilde\delta\delta )\left( *
  \left( \frac{g^{n-2k}}{(n-2k)!}R^{k-1}h\right)\right),
  \]
where $ (\delta\tilde\delta +\tilde\delta\delta )$ is the formal adjoint of the
Hessian type operator $(D\tilde D+\tilde D D)$.

Finally, using  Stoke's theorem we conclude that:
\begin{gather*}
H_{2k}'\cdot h=\int_M\left(h_{2k}'\cdot h+\frac{h_{2k}}{2}{\rm
tr}_gh\right)\mu_g =-\frac{1}{2}\langle\frac{c^{2k-1}}{(2k-1)!}R^k,h\rangle +\frac{h_{2k}}{2}\langle g,h\rangle \\
\phantom{H_{2k}'\cdot h} =\frac{1}{2}\langle h_{2k}g-\frac{c^{2k-1}}{(2k-1)!}R^k,h\rangle
=\frac{1}{2}\langle T_{2k},h\rangle.\tag*{\qed}
\end{gather*}\renewcommand{\qed}{}
\end{proof}

\section{Applications}

\subsection[A generalized Yamabe problem \cite{Labbivariation}]{A generalized Yamabe problem \cite{Labbivariation}}

It results from the previous theorem (see \cite{Labbivariation}) that for a compact Riemannian $n$-manifold $(M,g)$ with $n>2k$, the  Gauss--Bonnet curvature
  $h_{2k}$ is constant if and only if the metric $g$ is a critical point of the functional
  $H_{2k}$ when
restricted to the set ${\rm Conf}_0(g)$ of metrics pointwise conformal to $g$
and having the same total volume.

The previous result makes the following Yamabe-type problem plausible:
{\it In  each conformal class of
a fixed Riemannian metric on a smooth compact manifold with dimension $n>2k$
 there exists a metric with $h_{2k}$ constant.}

The previous problem is closely related to the recent $\sigma_k$-Yamabe problem of Viaclovsky
involving the symmetric functions of the Schouten tensor, see
\cite{Viaclovsky, Trudinger} and the references therein.

\subsection[Generalized Einstein manifolds \cite{Patterson, Labbivariation, LabbipqEinstein}]{Generalized Einstein manifolds \cite{Patterson, Labbivariation, LabbipqEinstein}}

Einstein metrics
are  the critical metrics of the total scalar curvature functional once restricted
to metrics of unit volume. Equivalently, a metric $g$ is Einstein if its  Ricci tensor $cR$ is proportional to $g$: $cR=\lambda g$.

In a similar way, the critical metrics of the total Gauss--Bonnet
curvature functional $H_{2k}$ once restricted to metrics with unit volume shall be called
$(2k)$-Einstein metrics. The $2$-Einstein metrics are nothing but the usual Einstein metrics.

 These critical metrics were studied f\/irst by Patterson \cite{Patterson} who proved that a locally irreducible symmetric
Riemannian metric is $(2k)$-Einstein for any positive integer $k$, and that harmonic Riemannian metrics are $2$ and $4$-Einstein.

The previous $(2k)$-Einstein metrics are  characterized by the condition that the  contraction of order $(2k-1)$
of Thorpe's
tensor $R^k$ is proportional to the metric, that is
\begin{gather*}
c^{2k-1}R^k=\lambda g.
\end{gather*}
More generally, for  $0<p<2q<n$, we shall say that a  Riemannian $n$-manifold
 is $(p,q)$-Einstein~\cite{LabbipqEinstein} if the $p$-th  contraction
 of Thorpe's tensor $R^q$
is proportional to the metric
 $g^{2k-p}$, that is
\begin{gather*}
c^pR^q=\lambda g^{2q-p}.
\end{gather*}

We recover the usual Einstein manifolds for $p=q=1$ and the previous $(2q)$-Einstein condition
for $p=2q-1$. The $(p,q)$-Einstein metrics are all critical metrics for the total Gauss--Bonnet curvature functional $H_{2q}$.

For all $p\geq 1$, $(p,q)$-Einstein implies $(p+1,q)$-Einstein. In particular, the metrics with constant $q$-sectional curvature (that is the
sectional curvature of $R^q$ is constant) are $(p,q)$-Einstein for all $p$.

On the other hand, the $(p,q)$-Einstein condition neither implies nor is
 implied by the $(p,q+1)$-condition  as shown by the following examples.

  Let $M$ be a 3-dimensional non-Einstein Riemannian manifold and $T^k$ be the $k$-dimensional f\/lat torus, $k\geq 1$,
 then the Riemann curvature tensor $R$ of the Riemannian product $N=M\times T^k$ satisf\/ies $R^q=0$ for $q\geq 2$.
In particular $N$ is $(p,q)$-Einstein for all $p\geq 0$ and $q\geq 2$ but it
 is not $(1,1)$-Einstein.

On the other hand, let  $M$ be a 4-dimensional Ricci-f\/lat but not f\/lat manifold
(for example a~$K_3$ surface endowed with the Calabi--Yau metric), then the Riemannian product $N=M\times T^k$ is $(1,1)$-Einstein
but not $(q,2)$-Einstein for any $q$ with $0\leq q\leq 3$.

The $(2q)$-Einstein condition, or equivalently the $(2q-1,q)$-Einstein condition, seems
 to be so weak to imply any topological restrictions on the manifold. However, for lower values of $p$ we have the
following obstruction result:

\begin{theorem}[\cite{LabbipqEinstein}]
 Let  $k\geq 1$ and  $(M,g)$ be a   $(1,k)$-Einstein  manifold (i.e.\ $cR^q=\lambda g^{2q-1}$)
of dimension
   $n\geq 4k$.
Then  the Gauss--Bonnet curvature $h_{4k}$ of  $(M,g)$ is nonnegative. Furthermore,
$h_{4k}\equiv 0$ if and only if  $(M,g)$ is $k$-flat.

In particular, a compact  $(1,k)$-Einstein  manifold of dimension   $n= 4k$  has its
 Euler--Poincar\'e characteristic nonnegative. Furthermore,
it is zero if and only if the metric is  $k$-flat.
\end{theorem}
The previous theorem generalizes a similar result of Berger about usual four dimensional Einstein manifolds.

\subsection[$(2k)$-minimal submanifolds \cite{Labbi2kminimal}]{$\boldsymbol{(2k)}$-minimal submanifolds \cite{Labbi2kminimal}}

Let $(\tilde M,\tilde g)$ be an $(n+p)$-dimensional Riemannian manifold,
 and let $M$ be an
$n$-dimensional  submanifold of $\tilde M$.

We shall characterize those submanifolds (endowed with the induced metric)
 that are critical points of the total Gauss--Bonnet curvature function.

Let $F$ be a local variation of $M$, that is a smooth map
\[ F:M\times (-\epsilon, \epsilon ) \rightarrow \tilde M,\]
such that $F(x,0)=x$ for all $x\in M$ and with  compact support ${\rm supp}F$.

The implicit function theorem implies that there exists $\epsilon >0$ such that for all
$t$ with $|t|<\epsilon$, the map
$\phi_t=F(\cdot,t):M\rightarrow \tilde M$ is a dif\/feomorphism onto a submanifold~$M_t$ of~$\tilde M$.

Let $g_t=\phi_t^*(\tilde g)$. Note that $g_1=g$.

\begin{theorem}[\cite{Labbi2kminimal}]
 Let
$\xi={\frac{d}{dt}}_{|t=0}\phi_t$ denotes the variation vector
field relative to a local variation~$F$  of $M$ with compact support as above.
\begin{enumerate}\itemsep=0pt
\item
 If
$H_{2k}(t)=\int_M h_{2k}(g_t)\mu_{g_t}$ denotes the total $(2k)$-th Gauss--Bonnet curvature of $\phi_t(M)$,
then
\begin{gather*}
H_{2k}'(0)=\int_M h_{2k+1}(\xi^\bot)\mu_g,
\end{gather*}
where $h_{2k+1}$ is the $(2k+1)$-th Gauss--Bonnet curvature of $M$  defined by~\eqref{GaussBonnetodd}.

\item
The submanifold $M$ is a critical point for the total
$(2k)$-th Gauss--Bonnet
curvature function for all local variations of $M$ if and only if the $(2k+1)$-Gauss--Bonnet
curvature $h_{2k+1}(N)$ of $M$  vanishes for all normal directions $N$.
\end{enumerate}

\end{theorem}

With reference to the previous variational formula and by analogy to the case of
usual minimal submanifolds we set the following def\/inition:
\begin{definition} For $0\leq 2k\leq n$,
an $n$-submanifold $M$ of a Riemannian manifold $(\tilde M,\tilde g)$ is said to be $(2k)$-minimal if
$h_{2k+1}\equiv 0.$
\end{definition}

Note that since $h_{2k+1}(N)=\langle T_{2k},B_N\rangle $, a submanifold is $(2k)$-minimal if and only if $T_{2k}$ is orthogonal to  $B_N$ for all
normal directions $N$.
Note the analogy with usual minimal submanifolds ($T_0=g$).

We list below some examples:
\begin{enumerate}\itemsep=0pt
\item  A f\/lat submanifold is always $(2k)$-minimal for all $k> 0$.  In fact
$R\equiv 0\Rightarrow
h_{2k+1}\equiv 0$. This shows that $(2k)$-minimal does not imply the usual minimality
condition.
\item  A totally geodesic submanifold is always $(2k)$-minimal for all $k\geq 0$.  In fact $B\equiv 0\Rightarrow
h_{2k+1}\equiv 0$.
\item If $M$ is a hypersurface of the Euclidean space then $(2k)$-minimality coincides with
Reilly's $(2k)$-minimality, \cite{Reillyrminimal}. On the other hand, if $M$ is a hypersurface of a space form $(\tilde M,\tilde g)$ of constant $\lambda$ then
$M$ is $(2k)$-minimal if and only if
\[ \sum_{i=0}^k\frac{(2k-2i+1)!(n-2k-1+2i)!\lambda^i}{i!(k-i)!}s_{2k-2i+1}=0.\]
In particular, $M$ is 2-minimal if and only if $6s_3+(n-1)(n-2)s_1\lambda =0$.
 Notice the dif\/ference with Reilly's
 $r$-minimality.
\item
A complex submanifold $M$ of a Kahlerian manifold $(\tilde M,\tilde g)$ is
 $(2k)$-minimal for any $k$.
\end{enumerate}

 Let now $f$ be a smooth function on $(M,g)$. We def\/ine the $\ell_{2k}$-Laplacian \cite{Labbi2kminimal}
 operator
 of $(M,g)$ as
 \begin{gather}\label{l2kLaplacian}
 \ell_{2k}(f)=-\langle T_{2k},{\rm Hess}(f)\rangle,
 \end{gather}
where $T_{2k}$ denotes the $(2k)$-th Einstein--Lovelock tensor (\ref{Einstein-Lovelock}) of $(M,g)$ and
 $0\leq 2k<n$, ${\rm Hess}(f)$ is the Hessian of $f$.

For $k=0$ we have $T_0=g$ and then $\ell_0=\Delta$ is the usual Laplacian.

For a compact manifold, the generalized Laplacian $\ell_{2k}$ satisf\/ies the following interesting properties:

For each $k\geq 0$, $\ell_{2k}(f)$ is a divergence hence
$\int_M\ell_{2k}(f)dv\equiv 0$. Furthermore, the operator $\ell_{2k}$ is self adjoint with respect to
the integral scalar product.

If for some $k$ with $0\leq 2k<n$,
the Einstein--Lovelock tensor $T_{2k}$ is positive def\/inite (or negative def\/inite), then
the operator $\ell_{2k}$ is elliptic and positive def\/inite (resp. negative def\/inite).

We shall say that the function $f$ is $\ell_{2k}$-harmonic if $\ell_{2k}(f)=0.$ In \cite{Labbi2kminimal} we proved the following maximum principle:
\begin{theorem}[\cite{Labbi2kminimal}]
Let $(M,g)$ be a compact manifold of positive definite (or negative definite) Einstein--Lovelock
tensor $T_{2k}$ then every  smooth and $\ell_{2k}$-harmonic function on $M$ is constant.
\end{theorem}

As a consequence of the previous result we proved the following about $(2k)$-minimal submani\-folds of the Euclidean space:
\begin{theorem}
A submanifold $M$ of the Euclidean space is $(2k)$-minimal if and only if the coordinate
functions restricted to $M$ are $\ell_{2k}$-harmonic functions on $M$.
\end{theorem}

\begin{corollary}
Let $0\leq 2k<n$ and let $(M,g)$ be a compact Riemannian  $n$-manifold with
 posi\-tive definite (or negative definite) Einstein--Lovelock tensor $T_{2k}$. Then there is
no non trivial isometric $(2k)$-minimal immersion of $M$ into the Euclidean space.
\end{corollary}

Note that the condition of positive (or negative) def\/initeness of $T_{2k}$ in the previous corollary is
necessary, as the f\/lat torus admits (non trivial) $(2k)$-minimal isometric immersions into the Euclidean space.

\subsection*{Acknowledgments}
The author would like to thank the referees for useful comments and especially  for indicating me the related work of
Patterson \cite{Patterson}.

\pdfbookmark[1]{References}{ref}
\LastPageEnding


\begin{thebibliography}{99}

\footnotesize\itemsep=0pt

\bibitem{Adler} Adler R.J., On excursion sets, tube formulas and maxima of random f\/ields,
{\it Ann. Appl. Probab.} {\bf 10} (2000), 1--74.

\bibitem{Albin} Albin P., Renormalizing curvature integrals on Poincar\'e--Einstein manifolds,
\href{http://arxiv.org/abs/math.DG/0504161}{math.DG/0504161}.

\bibitem{Berger} Berger M., Quelques formules de variation pour une
structure riemannienne, {\it Ann. Sci. Ecole Norm. Sup.~(4)} {\bf 3} (1970),  285--294.

\bibitem{Berningvariation} Bernig A., Variations of curvatures of subanalytic spaces and Schl\"af\/li-type formulas, {\it Ann. Global Anal. Geom.} {\bf 24} (2003), 67--93.

\bibitem{Berningreview} Bernig A., On some aspects of curvature, available at \url{http://homeweb1.unifr.ch/BernigA/pub/}.

\bibitem{Besse} Besse A.L., Einstein manifolds, Springer-Verlag, Berlin, 1987.


\bibitem{Bourguignon}   Bourguignon J.P., Les vari\'et\'es de dimension $4$ \`a
signature non nulle dont la courbure est harmonique sont d'Einstein, {\it Invent. Math.} {\bf 63} (1981), 263--286.

\bibitem{Chijap} Cai R.G., Ohta  N., Black holes in pure lovelock gravities,
{\it Phys. Rev. D} {\bf74} (2006), 064001, 8~pages, \mbox{\href{http://arxiv.org/abs/hep-th/0604088}{hep-th/0604088}}.


\bibitem{Cheeger} Cheeger J., M\"{u}ller  W., Schr\"{a}der  R., On the curvature of piecewise
f\/lat spaces, {\it Comm. Math. Phys.} {\bf 92} (1984), 405--454.

\bibitem{Deruelle-Madore} Deruelle  N., Madore J., On the quasi-linearity of the Einstein--Gauss--Bonnet gravity f\/ield equations,
\mbox{\href{http://arxiv.org/abs/gr-qc/0305004}{gr-qc/0305004}}.

\bibitem{Fenchel} Fenchel W., On total curvatures of Riemannian manifolds.~I,
{\it J. Lond. Math. Soc.} {\bf 15} (1940), 15--22.

\bibitem{Gray} Gray A., Some relations between curvature and characteristic classes,
{\it Math. Ann.} {\bf 184} (1970), 257--267.

\bibitem{Graybook} Gray  A., Tubes,  {\it Progress in Mathematics},  Vol.~221, Birkh\"auser Verlag, Basel, 2004.

\bibitem{Japweyl} Ishihara T., Kinematic formulas for Weyl's curvature invariants on submanifolds in complex projective space, {\it Proc. Amer. Math. Soc.} {\bf 97} (1986), 483--487.


\bibitem{Kowalski} Kowalski  O., On the Gauss--Kronecker curvature tensors,
{\it Math. Ann.} {\bf 203}  (1973), 335--343.


\bibitem{Kulkarni} Kulkarni  R.S.,
On Bianchi identities, {\it Math. Ann.} {\bf  199} (1972), 175--204.


\bibitem{Labbidoubleforms} Labbi M.L.,  Double forms, curvature structures and the
$(p,q)$-curvatures, {\it Trans. Amer. Math. Soc.} {\bf 357} (2005),
 3971--3992, \href{http://arxiv.org/abs/math.DG/0404081}{math.DG/0404081}.

\bibitem{LabbiPacific} Labbi M.L., On compact manifolds with positive second Gauss--Bonnet curvature, {\it Pacific J. Math.} {\bf 227} (2006),  295--310.

\bibitem{Labbivariation} Labbi M.L.,  Variational properties of the
Gauss--Bonnet curvatures, in Calculus of Variations and Partial Dif\/ferential Equations,
to appear, \href{http://arxiv.org/abs/math.DG/0406548}{math.DG/0406548}.

\bibitem{LabbiHabilitation} Labbi M.L., Riemannian curvature: variations on dif\/ferent notions of positivity, \href{http://arxiv.org/abs/math.DG/0611371}{math.DG/0611371}.

\bibitem{LabbiWeitzenbock} Labbi M.L., On Weitzenb\"ock curvature operators,
\href{http://arxiv.org/abs/math.DG/0607521}{math.DG/0607521}.

\bibitem{LabbipqEinstein} Labbi M.L.,  Remarks on  generalized Einstein manifolds,
\href{http://arxiv.org/abs/math.DG/0703028}{math.DG/0703028}.


\bibitem{Labbi2kminimal} Labbi M.L., On $(2k)$-minimal submanifolds,
\href{http://arxiv.org/abs/0706.3092}{arXiv:0706.3092}.


\bibitem{Lafontaine} Lafontaine J., Mesures de courbure des vari\'et\'es
lisses et des poly\`edres [d'apr\`es Cheeger, M\"{u}ller et Schr\"{a}der], {\it S\'eminaire  Bourbaki,
38\`eme ann\'ee, 1985--86}, no.~664, 1986.



\bibitem{Lovelock} Lovelock D., The Einstein tensor and its
generalizations, {\it J. Math. Phys.} {\bf  12} (1971), 498--501.

\bibitem{Madore} Madore J., Cosmological applications of the Lanczos Lagrangian,
{\it Classical Quantum Gravity} {\bf 3} (1986), 361--371.

\bibitem{Muto} Muto Y., Critical Riemannian metrics, {\it Tensor} {\bf 29} (1975), 125--133.

\bibitem{Nasu} Nasu T., On conformal invariants of higher order,
{\it Hiroshima Math. J.} {\bf 5} (1975), 43--60.

\bibitem{Nomizu} Nomizu  K.,
On the decomposition of generalized curvature tensor f\/ields. Codazzi, Ricci, Bianchi
 and Weyl revisited,
in Dif\/ferential Geometry (in Honor of Kentaro Yano), Kinokuniya, Tokyo, 1972, 335--345.


\bibitem{Patterson} Patterson E.M., A class of critical Riemannian metrics, {\it J. London Math. Soc.} {\bf 2} (1981), 349--358.



\bibitem{Reillyrminimal} Reilly  R.C., Variational properties of functions of the mean curvatures
for hypersurfaces in space forms, {\it J.~Differential Geom.} {\bf 8} (1973), 465--477.

\bibitem{Reillyproc} Reilly R.C., Variational properties of mean curvatures, in J. Proc. Summer
Sem. Canad. Math. Congress, 1971, 102--114.


\bibitem{Senovilla}  Senovilla  J.M.M., Super-energy tensors,
{\it Classical Quantum Gravity} {\bf 17} (2000),  2799--2842, \href{http://arxiv.org/abs/gr-qc/9906087}{gr-qc/9906087}.



\bibitem{Trudinger} Sheng W., Trudinger  N.S., Wang  X.-J.,
The Yamabe problem for higher order curvatures, {\it J. Differential Geom.}, to appear.


\bibitem{Stehney} Stehney A., Courbure d'ordre $p$ et les classes de Pontrjagin,
{\it J. Differential Geom.} {\bf 8} (1973), 125--134.



\bibitem{Thorpe} Thorpe  J.A., Some remarks on the Gauss--Bonnet
integral, {\it J. Math. Mech.} {\bf 18} (1969), 779--786.

\bibitem{Viaclovsky} Viaclovsky  J., Conformal geometry and fully nonlinear equations,
in World Scientif\/ic Memorial Volume for S.S.~Chern, to appear, \href{http://arxiv.org/abs/math.DG/0609158}{math.DG/0609158}.



\bibitem{Weyl} Weyl H.,  On the volume of tubes, {\it Amer. J.
Math.} {\bf 61} (1939), 461--472.






\end{thebibliography}
\end{document}